\documentclass[11pt,leqno]{article}
\usepackage{graphicx, amsfonts, amsthm, amsxtra, amssymb, verbatim}
\usepackage[mathscr]{euscript}
\begin{document}

\begin{center}
{\LARGE\bf Painlev\'e VI equations with algebraic solutions and family of curves \footnote{ 
Math. classification:  34M55, 35Q53 
}}
\\
\vspace{.25in} {\large {\sc Hossein Movasati, Stefan Reiter}} \\
Instituto de Matem\'atica Pura e Aplicada, IMPA, \\
Estrada Dona Castorina, 110,\\
22460-320, Rio de Janeiro, RJ, Brazil, \\
E-mail:
{\tt hossein@impa.br, reiter@impa.br} 
\end{center}
\begin{abstract}
In  families of Painlev\'e VI differential equations having common algebraic solutions we classify all the members
which come from geometry, i.e. the corresponding linear differential equations which are Picard-Fuchs associated to families of
algebraic varieties. In our case, we have one family with zero dimensional fibers and all others are families of curves.
We use the classification of families of elliptic curves with four singular fibers done by Herfurtner in 1992 and generalize the results of
Doran in 2001 and Ben Hamed and Gavrilov in 2005.  
\end{abstract}
\newtheorem{theo}{Theorem}
\newtheorem{exam}{Example}
\newtheorem{coro}{Corollary}
\newtheorem{defi}{Definition}
\newtheorem{prob}{Problem}
\newtheorem{lemma}{Lemma}
\newtheorem{prop}{Proposition}
\newtheorem{rem}{Remark}
\newtheorem{conj}{Conjecture}
\def\Gal{{\rm Gal}}              
\def\Z{\mathbb{Z}}                   
\def\Q{\mathbb{Q}}                   
\def\C{\mathbb{C}}                   
\def\ring{{\sf R}}                             
\def\R{\mathbb{R}}                   
\def\N{\mathbb{N}}                   
\def\uhp{{\mathbb H}}                
\newcommand{\mat}[4]{
     \begin{pmatrix}
            #1 & #2 \\
            #3 & #4
       \end{pmatrix}
    }                                

\newcommand\SL[2]{{\rm SL}(#1, #2)}    
\def\per{{\sf pm}}

\def\la{\lambda}
\def\th{\theta}
\def\P{\mathbb P}
\theoremstyle{plain}
\def\NN{{\Bbb N}}
\def\CC{{\Bbb C}}
\def\ZZ{{\Bbb Z}}
\def\QQ{{\Bbb Q}}
\def\SL{{\rm SL}}
\def\GO{{\rm GO}}
\def\GL{{\rm GL}}
\def\PGL{{\rm PGL}}
\def\dz{{\rm d} z}
\def\dx{{\rm d} x}

\def\la{{\lambda}}
\def\nlambda{{\tilde{\lambda}}}
\def\nue{{\nu}}
\def\ele{{\rm L}}
\def\tparam{{t}}
\def\Mu{{\mu}}
\def\th{{\theta}}
\def\thh{{\tilde{\theta}}}
\def\Th{{\Theta}}
\def\tr{{\rm tr}}
\def\rk{{\rm rk}}
\def\Trace{{\rm Tr}}
\def\Mat{{\rm Mat}}
\def\al{{\alpha}}
\def\diag{{\rm diag}}
\def\id{{\rm id}}
\def\k{{\frak k}}
\def\l{{\frak  l}}

\section{Introduction}
Along the solutions of the  sixth Painlev\'e differential equation written in the vector field form
\begin{equation}
\label{paiham}
PVI_\th:\ \frac{\partial K}{\partial \mu} \frac{\partial}{\partial \lambda}-
\frac{\partial K}{\partial \la}\frac{\partial}{\partial \mu}+
\frac{\partial}{\partial t}
\end{equation}
in $\C^3$ with coordinates $(\la,\mu,t)$, where
$$
t(t-1)K= \la (\la-1)(\la-t)\mu^2-(\th_2(\la-1)(\la-t)+\th_3 \la(\la-t)+
   (\th_1-1)\la(\la-1))\mu+\kappa (\la-t),
$$
$$
\kappa=\frac{1}{4}((\sum_{i=1}^3\th_i-1)^2-\th_4^2),
$$
and $\th=(\th_1,\th_2,\th_3,\th_4)$ is a fixed multi-parameter, the linear differential equation 
\begin{equation}
 \label{11.4.08}
y''+p_1(z)y'+p_2(z)y=0
\end{equation}
$$
p_1(z):=\frac{1-\th_1}{z-t}+\frac{1-\th_2}{z}+\frac{1-\th_3}{z-1}-\frac{1}{z-\la},\
$$
$$
p_2(z):=\frac{\kappa}{z(z-1)}-
\frac{t(t-1)K}{z(z-1)(z-t)}+\frac{\la(\la-1)\mu}{z(z-1)(z-\la)}
$$
is isomonodromic, i.e. its monodromy group representation is constant.
Recently, there have been works on algebraic solutions of
(\ref{paiham}) using linear differential  equations coming from geometry, i.e. those who  are Picard-Fuchs equations associated 
to families of varieties. Linear equations (\ref{11.4.08}) with finite monodromy come automatically from geometry and this is the origin of many algebraic
solutions known until now (see \cite{boa06} and the references therein). Doran in 
\cite{dor01} took $5$ deformable families of elliptic curves with non-constant $j$ invariant along the deformation parameter and with
exactly four singular fibers, which were already classified by Herfurtner in \cite{her91}, and obtained
algebraic solutions for (\ref{paiham}) (see Table 1, column 1 and 2, row 3-6 for $a=c=\frac{1}{2}$). 
In general one can obtain such algebraic solutions by taking
pull-backs of the Gauss hypergeometric equation (see \cite{kit05} and the references therein).  
Ben Hamed and Gavrilov in \cite{ben05} took zero dimensional families of three points
varieties, constructed directly from the Herfurtner list, and they obtained the same algebraic solutions in 
the $(t, \lambda)$-space but for different parameters $\th$. 
They have also noticed that  in the parameter space  of Painlev\'e equations (\ref{paiham}), the points obtained by them and Doran lie in families with algebraic solutions whose projections in the $(t,\lambda)$ space is independent of the parameter of the family. 
Then by a straightforward calculation they showed that up to
the Okamoto transformations corresponding to the M\"obius transformation of $\P^1$, such families of Painlev\'e equations 
are given by the first and second column of Table 1. 
The main result of this article is to classify the members of such families which come from geometry (third column of Table 1). 
{\small
\begin{center}
Table 1: Algebraic solutions of families of the sixth Painlev\'e equation
\begin{tabular}{|c|c|c|}
\hline
Algebraic solution  &   $(\th_1,\th_2,\th_3,\th_4)$  &   family of algebraic varieties  \\ \hline
\hline
 $\{\lambda=t=0\}\cup \{\lambda=t=1\} $  &   $(0,1-c,$ &  $y=x^{1-a}(1-x)^b(z-x)^{1+a-c},\ \frac{dx}{y}$\\  & $c-a-b,$ &     \\ & $b-a)$ &  \\\hline
  $\la=(\frac{-a+1}{a+2c-3})b,\ t=b^2$ &   $(\frac{1}{2},a-1,$ & zero dimensional varieties\\ $\mu=\frac{-a-2c+3}{2b}$ & $\frac{1}{2},-(a+2c-3))$ &   \\ \hline
   $\la=-b, \ t=b^2$ &   $(c-\frac{1}{2}, a+c-1,$ & $y=(4x^2-g_2x+g_3)^c(x+g_2/4)^a,\ \frac{dx}{y}$  \\ $\mu=\frac{-a-2c+2}{2b}$& $ c-\frac{1}{2}, a+c-1)$  &  $g_2=4(\tilde z^2+\tilde z)$
 \\ 
& &  $g_3=-12\tilde b\tilde z^3-8\tilde z^4+8\tilde z^3-8\tilde z^2$\\ 
 & & $\tilde b=\frac{3}{4}(b+\frac{1}{b})+\frac{1}{2} , \tilde z=\frac{-1}{b}z$\\ \hline
 $\la=\frac{-2b-1}{b^2}, \ t=\frac{2b+1}{b^4+2b^3}$ &$(a-\frac{1}{2}, 3(a-\frac{1}{2}), $ &  $y=(4x^3-g_2x-g_3)^a,\ \frac{dx}{y}$\\ & $a-\frac{1}{2}, a-\frac{1}{2})$ &  $g_2=12\tilde z ^2(\tilde z ^2+2\tilde b \tilde z +1)$ \\ 
$\mu=\frac{(-2a+3)b^2(b+2)}{2(b+1)^2}$ &    &   $g_3=4\tilde z ^3(2\tilde z ^3+3(\tilde b ^2+1)\tilde z ^2+6\tilde b \tilde z +2)$ \\
& &  $\tilde b =\frac{2}{3}(b+\frac{1}{b})-\frac{1}{3}$, $\tilde z= -\frac{b^2+2b}{3} z$  \\ \hline
 $\la=\frac{b^3+b^2+3b+3}{b^3+b^2-5b+3},\ t=\frac{b^4-6b^2-8b-3}{b^4-6b^2+8b-3}$ & $(a-\frac{1}{2}, \frac{1}{2},$ & $y=(4x^3-g_2x-g_3)^a,\ \frac{dx}{y}$ \\ & $ a-\frac{1}{2}, a-\frac{1}{2})$&  $g_2=3\tilde z^3(\tilde z+\tilde b),\ g_3=\tilde z^5(\tilde z+1)$\\
$\mu=\frac{(3a-2)(b-1)^2(b+3)}{24(b+1)}$ &     &  $\tilde b=\frac{2}{3}\frac{b^2-3}{b^2+3}+\frac{1}{3}, \tilde z=-\frac{b^3-3b^2+3b-1}{b^3-3b^2+3b-9}z$  
 \\ \hline
 $\la=\frac{-2b^2-4}{b^4-6b^2},\
t=\frac{-12b^2+8}{b^6-6b^4}$  & $(a-\frac{1}{2}, \frac{1}{3},$ & $y=(4x^3-g_2x-g_3)^a,\ \frac{dx}{y}$\\ & $a-\frac{1}{2}, 2a-1)$ &  $g_2=3\tilde z^3(\tilde z+2\tilde b)$\\ $\mu=\frac{(-3a+2)b^2(b^2+2)(b^2-6)}{12(b^2-2)^2}$ 
& &  $g_3=\tilde z^4(\tilde z^2+3\tilde b\tilde z+1)$ \\
 &     & $\tilde b=\frac{1}{4}(b+\frac{2}{b}) ,\ \tilde z=-\frac{2b^3}{3b^2-2}z $
\\ \hline
\end{tabular}
\end{center}
}

\begin{theo}
In  the columns 1 and 2 of Table 1 the corresponding linear differential equations (\ref{11.4.08}) come from geometry if and only the exponent 
parameters $a$ and $c$ 
in column 2 (and $b$ in the first row)  are rational numbers. The corresponding family of algebraic varieties and differential form are 
listed in column 3. Moreover,
the second family (together with its algebraic solution) is Okamoto equivalent to the third family and the fourth family is Okamoto equivalent to the fifth family. Both Okamoto transformations are equivalent to the middle convolution of the corresponding Fuchsian systems.
\end{theo}
In the first row of Table 1, the corresponding linear differential equation (\ref{11.4.08}) is the 
Gauss hypergeometric equation and  the geometric interpretation is classical  and it has nice applications in the theory of the special values of Gauss hypergeometric functions (see \cite{shi04}). Note that in this case the projection of the
corresponding  algebraic curve in the $(\lambda,\mu,t)$ space into the $(\lambda,t)$ space is just the zero dimensional variety $\{(0,0),(1,1)\}$. In all other cases it is a curve in the $(\lambda,t)$ space. In the literature one finds mainly 
the equations of such curves.  
Note also that the parametrization of the algebraic solutions in column 1 is different from the 
parametrization of the solutions of the vector field  (\ref{paiham}). 

Let us put the parameters $a$ and $c$  in column $2$ and rows $3,4,5,6$ equal to $\frac{1}{2}$.  We 
obtain families of elliptic curves $y^2=4x^3-g_2(z,b)x-g_3(z,b)$ with exactly four singular fibers and with $j$-invariant
depending on an extra parameter $b$. These are exactly four families of the five Herfurtner families. The missing
family in the Herfurtner list is the one given by 
$
g_2=3(z-1)(z-b^2)^3, \ g_3=(z-1)(z-b^2)^4(z+b).
$ 
The corresponding family of curves
$y=(4x^3-g_2x-g_3)^{\tilde a}$ gives us the Painlev\'e equation and its algebraic solution in Table 1, row 3 by setting
$c=\frac{5}{6},\ a=\tilde a-\frac{1}{3}$. 
Note that in this case  we find two apparently different geometric interpretations 
for the same Painlev\'e equation. Note also that in the mentioned five families of elliptic curves just for the family used in row $3$ the polynomial  $4x^3-g_2(z,b)x-g_3(z,b)$ is reducible in $x$.  
In row $2$ for $a$ and $c$ rational numbers we have shown that the monodromy group of the linear equation (\ref{11.4.08})
is a dihedral group and so it is finite. Also the other families are related via the middle convolution to (third order) differential equations whose monodromy groups are
finite imprimitive reflection groups for rational parameters. In \cite{Boalch03} Boalch started with third order differential equations to obtain (via a construction equivalent to the middle convolution) algebraic solutions and the
parameters of the corresponding Painlev\'e VI differential equations.   
Recently Cantat and Loray showed in \cite[Prop. 5.4]{CL07} that any algebraic solution of a Painlev\'e VI
differential equation having degree 2, 3, or 4 belongs (up to Okamoto transformation) to one of the families in Table 1.

Using Singular (see \cite{GPS01}), for each linear equation (\ref{11.4.08}) we have calculated the corresponding Fuchsian system in the Schlesinger form and the Okamoto transformation of the algebraic solutions in Table 1 corresponding to exchanging the role of the first and second coordinates of the system. The details of the calculation is explained for the example in 
row 3 of Table 1 and for the others the reader is referred to the home-page of the first author.

Let us explain the content of each section. In \S \ref{t2t3} we introduce systems of linear differential equations in 
two and three variables. Pulling these back  we get Fuchsian systems with four singularities.
In \S \ref{review} we recall some well-known facts about linear differential equations and in \S \ref{algorithm} we explain
how to calculate Fuchsian systems in Schlesinger form associated to families of curves in Table 1, column 3. Such calculations
are explained for the third row of Table 1 in \S \ref{example}. 
In \S\ref{middleconvolution} we recall some basic facts about the middle convolution and in \S \ref{midconvI} we
show that the algebraic solution in  row 2 of Table 1 can be obtained via the middle convolution of 
the Fuchsian system corresponding to the algebraic solution in row 3 of Table 1. In \S\ref{monodromy} we show that
the Fuchsian system associated to the Painlev\'e equation in row 2 of Table 1 has finite monodromy if $a,c\in \Q$ and
so it comes from variation of zero dimensional varieties. Further we point out how all the other Picard-Fuchs differential equations are related via the middle convolution to those with finite monodromy.
Finally, in \S \ref{midconvII} 
we show that the middle convolution relates the results in Table 1 row 4 and row 5.\\

The second author gratefully acknowledges financial support by the CNPq
and thanks IMPA for its hospitality. The authors thank L. Gavrilov for useful discussions in the initial
steps of the present article.

\section{Linear systems in two variables}
\label{t2t3}
For $a,b,c\in\C$ fixed, we consider the following family of  transcendent curves:
\begin{equation}
\label{khodaya1}
E: y=f(x), \ 
\end{equation}
$$
f(x):=(t_1-t_3)^{\frac{1}{2}(1-a-c)}(t_1-t_2)^{\frac{1}{2}(1-a-b)}(t_2-t_3)^{\frac{1}{2}(1-b-c)}
(x-t_1)^a(x-t_2)^b(x-t_3)^c 
$$
Here  $t=(t_1,t_2,t_3)$ is a parameter in
$$
T:=\{t\in \C^3\mid (t_1-t_2)(t_2-t_3)(t_3-t_1)\not =0\}.
$$
 We distinguish three, not necessarily closed,  paths  in $E$.
In the $x$-plane let $\tilde \delta_i, \ i=1,2,3,$ be the straight path connecting
$t_{i+1}$ to $t_{i-1}$, $i=1,2,3$ (by definition $t_4:=t_1$ and $t_0=t_3$). 
There are many paths in $E$
which are mapped to $\tilde \delta_i$ under the projection on $x$. 
We choose one of them
and call it $\delta_i$. 
We can make our choices so that $\delta_1+\delta_2+\delta_3$ is a limit of a closed and 
homotopic-to-zero path in $E$. For instance, we can take the path $\tilde \delta_i$'s in
such a way that the triangle   formed by them has almost zero area.  Now, we have the integral
\begin{equation}
\label{4nov06}
\int_{\delta}\frac{p(x)dx}{y}=\int_{\tilde \delta}\frac{p(x)dx}{f(x)},\ p\in\C[x],
\end{equation}
where $\delta$ is one of the paths explained above.
By a linear change in the variable $x$ such integrals can be written in terms of the Gauss hypergeometric function (see \cite{iwa91}).
Another way to study the integrals (\ref{4nov06}) is by using Pochhammer cycles. 
For simplicity
we explain it for the pairs $(t_1,t_2)$. 
The Pochhammer cycle associated to the points $t_1,t_2\in\C$ and the path 
$\tilde \delta_3$ is the commutator
$$
\tilde\alpha_3=[\gamma_{1},\gamma_{2}]=\gamma_{1}^{-1}\cdot \gamma_{2}^{-1}\cdot \gamma_{1}\cdot 
\gamma_{2},
$$
where $\gamma_{1}$ is a loop along $\tilde \delta_3$ starting and ending at some
point in the middle of $\tilde \delta_1$ which encircles $t_{1}$ once anti-clockwise, 
and $\gamma_{2}$ is a similar  loop with respect to $t_{2}$. 
It is easy to see that the cycle $\tilde \alpha_3$ lifts up to a closed
path $\alpha_3$ in $E_t$ and
if $a,b\not \in \Z$
then 
$$\int_{\alpha_3}\frac{p(x)dx}{y}=
( 1-e^{-2\pi i a})
( 1-e^{-2\pi i b})
\int_{\tilde \alpha_3} \frac{p(x)}{f(x)}dx.
$$ 
(see \cite{iwa91}, Proposition 3.3.7). 

For a fixed ${\sf a}\in T$, the period map is given by:
\begin{equation}
\label{permap}
\per: (T,{\sf a})\rightarrow \GL(2,\C),\ t\mapsto
\mat
{\int_{\delta_1}\frac{dx}{y}}
{\int_{\delta_2}\frac{dx}{y}}
{\int_{\delta_1}\frac{xdx}{y}}
{\int_{\delta_2}\frac{xdx}{y}},
\end{equation}
where $(T,{\sf a})$ means a small neighborhood of ${\sf a}$ in $T$. 
The map $\per$ can be extended along any path in $T$ with the starting point 
${\sf a}$. 
The period map $\per$ is a fundamental  system for the linear differential equation $dY=AY$ 
in $\C^3$, where
\begin{equation}
\label{A=}
A=
\end{equation}
{\tiny
$$
\frac{1}{(t_1-t_2)(t_1-t_3)}
\mat
{
\frac{1}{2}(b+c-2)t_1+ \frac{1}{2}(a+c-1)t_2+\frac{1}{2}(a+b-1)t_3 
}
{
-a-b-c+2
}
{
at_2t_3+(b-1)t_1t_3+(c-1)t_1t_2
}
{
-\frac{1}{2}(b+c-2)t_1-\frac{1}{2}(a+c-1)t_2
-\frac{1}{2}(a+b-1)t_3
}dt_1
$$
}
$$
+(\cdots)dt_2+(\cdots)dt_3
$$
The matrix coefficient of $dt_2$ (resp $dt_3$) is obtained by permutation of $t_1$ with $t_2$ 
and $a$ with $b$ (resp. $t_1$ with $t_3$ and $a$ with $c$) in the matrix coefficient of $dt_1$ written above.
Now, for the multi-valued function
$$
y=
(27t_3^2-t_2^3)^{\frac{1}{2}(\frac{1}{2}-a)}(4x^3-t_2x-t_3)^{a}
$$
we have the system
\begin{equation}
 \label{t1t2one}
A=\frac{1}{27t_3^2-t_2^3}\left (
\begin{pmatrix}
\frac{1}{4}t_2^2 & -27at_3+18t_3\\
-\frac{9}{4}at_2t_3+\frac{3}{4} t_2t_3 & -\frac{1}{4} t_2^2
\end{pmatrix}
dt_2+\right.
\end{equation}
$$
\left.
\begin{pmatrix}
-\frac{9}{2}t_3 & 18at_2-12t_2\\
\frac{3}{2}at_2^2-\frac{1}{2}t_2^2 & \frac{9}{2}t_3
\end{pmatrix}
dt_3
\right )
$$
and for
$$
y= 
(t_2^2+2t_3)^{\frac{1}{2}(1-a-c)}(t_2^2-16t_3)^{\frac{1}{2}(\frac{1}{2}-c)}
(4x^2-t_2x+t_3)^{c}(x+\frac{1}{4}t_2)^a
$$
we have
\begin{equation}
\label{t1t2two}
A=\frac{1}{(t_2^2-16t_3)(t_2^2+2t_3)}
\end{equation}
$$
\left(
\mat
{
(6at_2t_3-6ct_2t_3-\frac{1}{2}t_2^3+5t_2t_3)
}
{
(-48at_3-96ct_3+96t_3)
}
{
(12at_3^2-3ct_2^2t_3+t_2^2t_3-4t_3^2)
}
{
(-6at_2t_3+6ct_2t_3+\frac{1}{2}t_2^3-5t_2t_3)
}dt_2+\right.
$$
$$\left.
\mat
{
(-3at_2^2+3ct_2^2+t_2^2+8t_3)
}
{
(24at_2+48ct_2-48t_2)
}
{
(-6at_2t_3+\frac{3}{2}ct_2^3-\frac{1}{2}t_2^3+2t_2t_3)
}
{
(3at_2^2-3ct_2^2-t_2^2-8t_3)
}dt_3\right)
$$
The calculation of the matrix $A$ for the elliptic case $a=b=c=\frac{1}{2}$ is classical and goes back to 
Griffiths \cite{gri66}. In fact the matrix $A$ in (\ref{A=}) can be also calculated from the well-known Fuchsian
system for hypergeometric functions, i.e. for the case $t_1=0,\ t_2=1, \ t_3=t$. The algorithms for calculating
such a matrix are explained in \cite{mov08} and the implementation of such algorithms in a computer can be found
in the first author's home-page.

\section{Review of Fuchsian differential equations}
\label{review}

Here we collect some basic facts about Fuchsian differential equations
which can all be found in \cite{iwa91}.
Let $D=\frac{d}{\dz}$ and
\begin{eqnarray}\label{d.e.} D^2 y+p_1(z) Dy+p_2(z)y=0\end{eqnarray} 
be a Fuchsian differential equation with regular singularities
at $t_1,\ldots,t_m \in \CC$ and $t_{m+1} =\infty.$
Then by  \cite[Prop. 4.2]{iwa91}
$$
p_1(z)=\sum_{i=1}^m \frac{a_i}{z-t_i},\ 
p_2(z)=\sum_{i=1}^m \frac{b_i}{(z-t_i)^2}+\sum_{i=1}^m \frac{c_i}{(z-t_i)},\;
a_i,b_i,c_i \in \CC
$$
where $p_2(z) \prod_{i=1}^m (z-t_i)^2$
is a polynomial in $\CC[z]$ of degree at most $2(m-1),$ i.e.
$\sum c_i=0.$
The following table of the regular singular points 
$t_1,\ldots,t_{m+1}$ and the exponents $s_i^1,s_i^2$ at $t_i$ is called
the Riemann scheme of the differential equation 
\begin{eqnarray*} \left( \begin{array}{cccccc}
       t_1 & \ldots & t_{m+1} \\
        s_1^1 & \ldots & s_{m+1}^1 \\
           s_1^2 & \ldots & s_{m+1}^2
        \end{array}\right),\end{eqnarray*}
where the sum of all exponents satisfies the Fuchs relation
$\sum_{i=1}^{m+1} \sum_{j=1}^2 s_j^i =m-1.$ 
The exponents at the singularity $t_i,\;i=1,\ldots,m,$ are the roots of $s(s-1)+a_i  s+b_i=0$
(see  \cite[p. 170]{iwa91}) and at $t_{m+1}$ the exponents satisfy 
\begin{eqnarray*} s(s+1)-(\sum_{i=1}^m a_i)s+(\sum_{i=1}^m b_i+(\sum_{i=1}^m c_i t_i))=0.\end{eqnarray*} 
The second order Fuchsian differential equation (\ref{d.e.})
can be transformed into $\SL$-form by
substituting $y$ by $f y,$ where $0\neq f$
 satisfies  $Df=-\frac{1}{2} p_1(z) f.$ 
Then by  \cite[p. 166]{iwa91} 
 \begin{eqnarray} \label{SL} 
D^2 y =p(z)y,& p(z)=-p_2(z)+\frac{1}{4}p_1(z)^2+\frac{1}{2} Dp_1(z).\end{eqnarray}
>From a two dimensional system $DY=QY,\quad Q=(q_{ij})$
(not necessarily a Fuchsian system), we obtain the second order differential equation (\ref{d.e.}) 
for the first coordinate $y_1$ of $Y=(y_1,y_2)^\tr$ (see  \cite[Lemma 6.1.1]{iwa91}) with
\begin{equation}
 \label{25apr}
p_1(z)= -D \log (q_{12}(z)) -\Trace (Q), \
                   p_2(z)= \det(Q(z))-D q_{11}+q_{11} D \log q_{12}.
\end{equation}
If $\la$ is a zero of $q_{12}$ of order $r$ and $\la \not \in \{t_1,\ldots,t_{m}\}$
then $z=\la$ is an apparent singular point with the exponents $0$ and $r+1$
(see  \cite[Lemma 6.1.2]{iwa91}).

\section{The algorithm}
\label{algorithm}

In this section we explain how to obtain the algebraic solutions
of the Painlev{\'e} VI differential equation starting from the
families of curves in Table 1 column 3, row 3, 4, 5, 6, which are constructed directly  from the
Herfurtner list of families of elliptic curves (see \cite{her91}). For the convenience of the reader 
we have listed the five families which we need:

\begin{center}
{\small Table 2: List of deformable families of elliptic curves $y^2=4x^3-g_2x-g_3$ with four \\
singular fibers and with non constant 
$j$ invariant along the deformation} 
\begin{tabular}{||c|c||}
\hline\hline
name & deformation \\ \hline\hline
1 & $g_2=3(z-1)(z-b^2)^3, \ g_3=(z-1)(z-b^2)^4(z+b)$\\ \hline
2 & $g_2=12z^2(z^2+bz+1), \ g_3=4z^3(2z^3+3bz^2+3bz+2)$\\ \hline
3 & $g_2=12z^2(z^2+2bz+1),\ g_3=4z^3(2z^3+3(b^2+1)z^2+6bz+2)$\\ \hline
4 & $g_2=3z^3(z+b),\ g_3=z^5(z+1)$\\ \hline
5 & $g_2=3z^3(z+2b),\ g_3=z^4(z^2+3bz+1)$\\\hline
\end{tabular}
\end{center}
Let us consider $f(x)=4x^3-g_2x-g_3$, 
where $g_2, g_3$ correspond to one of the five families of elliptic curves with parameter 
$b$ in the Herfurtner list.
Since the roots of the discriminant $\Delta=27g_2^3-g_3^2$ of $f$, 
which will be the singular points of the corresponding differential equation
for $\int_\delta \frac{\dx}{y}$, are not all rational functions in $b$ we
substitute $b$ by a suitable rational function in $b$ such that the transformed roots are rational in $b$. 
Such substitutions are done by the equalities $\tilde b=\cdots$ in Table 1 column 3.   
In the next step we check whether the polynomial $f(x)$ factorizes over
$\QQ(b,z)$.
 It turns out that this only happens for second family. 
In that case we have
\begin{equation}
 \label{factorize}
f(x)=(4x^2-g_2x+g_3)(x+\frac{g_2}{4}).\ 
\end{equation}
Substituting $g_2$ and $g_3$ in the list 2 (resp. 1, 3, 4, 5) of Table 2 for $t_2$ and $t_3$ in (\ref{t1t2two})
(resp.  (\ref{t1t2one})) we
obtain a  system $ DY=AY$. 
We can now compute the second order differential equation satisfied by 
 $\int_\delta \frac{\dx}{y},$ where $\delta$ is a Pochhammer cycle, using the formula \eqref{25apr}
and \eqref{SL}.
It always turns out that 
 \eqref{d.e.} has four 
singularities at $t_1,t_2=0,t_3$ and $t_4=\infty$ and 1 apparent singularity at
$\la$ with exponents $0$ and $2.$ 
Hence, the $\SL$-form can be written as follows  (\cite[p. 173]{iwa91})
\begin{eqnarray*} p(z)=\sum_{i=1}^3 \frac{a_i}{(z-t_i)^2} +\frac{a_4}{z(z-t_3)}+
               \frac{t_1(t_1-t_3)/t_3 \cdot L}{z(z-t_1)(z-t_3)}
                +\frac{3}{4} \frac{1}{(z-\la)^2}-
                \frac{\la (\la-t_3)/t_3 \cdot \nu}{z(z-t_3)(z-\la)},\end{eqnarray*}
where
\begin{eqnarray*} a_i=\frac{1}{4} (\th_i^2-1),\; i=1,2,3,&a_4=-\frac{1}{4}(\sum_{i=1}^3\th_i^2-\th_4^2-1)-\frac{1}{2},\\
L= t_3 ((p(z)-\frac{a_1}{(z-t_1)^2})(z-t_1))_{z=t_1},& 
\nu=-t_3 ((p(z)-\frac{3}{4(z-\lambda)^2})(z-\lambda))_{z=\lambda}.
\end{eqnarray*}
The corresponding Riemann-scheme is
\begin{eqnarray*} \left( \begin{array}{cccccc}
       t_1 & t_2 & t_3 &  \infty & \la\\
       \frac{1-\th_1}{2} &  \frac{1-\th_2}{2} & \frac{1-\th_3}{2} & \frac{-1-\th_4}{2} & -\frac{1}{2}\\
        \frac{1+\th_1}{2}&  \frac{1+\th_2}{2} & \frac{1+\th_3}{2} & \frac{-1+\th_4}{2} & \frac{3}{2}\\
        \end{array}\right).\end{eqnarray*}
The system 
 \begin{eqnarray}\label{Schles} DY=QY,& Q=\sum_{k=1}^3 \frac{Q_k}{z-t_k},\ Q_k=Q_k(t)=(q_{ij}^k)\in \Mat_n(\C),\ k=1,2,3
\end{eqnarray}
 of Schlesinger type with the 4 regular singularities
 at $t_1,t_2=0,t_3,t_4=\infty$ can now be determined: 
 We can assume that $\th_i$ and $0$ are the eigenvalues
 of $Q_i$ and that (if $\th_4 \neq 1$)  
 \begin{eqnarray*}  -\sum_{i=1}^3 Q_i(t)= \left( \begin{array}{cc}
                                    \al & 0 \\
                                         0 &\al+\th_4-1
                                  \end{array}\right),& 
 2\al+\sum_{i=1}^4 \th_i = 1.
  \end{eqnarray*}
Further  we  normalize the singularities via the M\"obius transformation $z\mapsto z \cdot t_3$ to $t=\frac{t_1}{t_3}, 0, 1$ and $\infty$. (Note that also the
apparent singularity is transformed to $\frac{\lambda}{t_3}$ which we will also denote again by $\la$.)
 Then by \cite[Prop. 6.3.1]{iwa91} the matrices $Q_i,i=1,2,3,$ can be expressed
 as follows:

 \begin{eqnarray}\label{Qi} Q_i= \left( \begin{array}{cc}
                                  M_i(W_i-W) & -M_i \\
                           -(W_i-W)(M_i(W-W_i)+\th_i) & \th_i-M_i(W_i-W)
                                  \end{array}\right),\ i=1,\ldots,3,\end{eqnarray}
\begin{eqnarray*}
\begin{array}{c|c|c}
      &M_i & W_i \\
&&\\
\hline
&&\\
  i=1 & \frac{\la-t}{t (t-1)}& \la (\la-1)(\mu+\frac{\alpha}{\la})\\
&&\\
\hline
&& \\
 i=2&\frac{\la}{t}&(\la-t)(\la-1)(\mu+ \frac{\alpha}{\la})-\frac{t \alpha}{\la}\\
&&\\
\hline
&&\\
i=3&\frac{\la-1}{1-t}&\la (\la-t)(\mu+\frac{\alpha}{\la})
\end{array} &
&\begin{array}{l}
\\
\\
\\
 (\th_4-1)W=\sum_{i=1}^3 W_i(M_iW_i-\th_i), \\
\\
\\
 \mu=\nu-\frac{1}{2}\sum_{i=1}^3 \frac{1-\th_i}{\la-t_i}, \\
\\
\\
 \al=-\frac{1}{2}(\sum_{i=1}^4 \th_i - 1).
\end{array}
\end{eqnarray*}
Note that a zero of $Q_{1,2}$ has order $r=1,$ which is the apparent
singularity $\la$  for the first coordinate.
If we compute the differential equation for the second coordinate or equivalently
if we  compute the differential equation for the first coordinate of the transformed system
\begin{eqnarray*}
 DY=TQT^{-1} Y,& T=\left(\begin{array}{cc}
                             0& 1\\
                             1 & 0
                   \end{array}\right),
\end{eqnarray*}
then $(\th_1,\ldots,\th_4)$ is transformed to $(\th_1,\ldots,\th_3,2-\th_4)$.
\section{An example}
\label{example}
We demonstrate how the algorithm in \S \ref{algorithm} yields the results for the 
third row of of Table 1.
We start with the second family in Table 2. 
The roots of the discriminant of $f$ are
\begin{eqnarray*} 0,\omega_1,\omega_2, &\omega_{1,2}=-\frac{1}{3}(2b-1\pm 2\sqrt{b^2-b-2})\end{eqnarray*}
Solving the diophantic equation $\alpha^2-\alpha-2=\beta^2$, 
we substitute $b$ by $\frac{3}{4}(b+\frac{1}{b})+\frac{1}{2}$ 
in $g_2$ and $g_3.$
Thus the new roots are
\begin{eqnarray*} t_1=-b, & t_2=0, & t_3=-\frac{1}{b}.\end{eqnarray*}
Since $f(x)=4x^3-g_2 x-g_3$ factorizes we get (\ref{factorize}) with 
the new coefficients
\begin{eqnarray*} g_2:=4 (z^2+z), &
 g_3:=\frac{-9 b^2 z^3-8 b z^4+2 b z^3-8 b z^2-9 z^3}{b},\end{eqnarray*}
Substituting $g_2$ and $g_3$ for $t_2$ and $t_3$ in (\ref{t1t2two}) we
obtain a system which does not fit in our paper and we do not write it here.
Computing the \SL-form \eqref{SL} of the differential equation
for the first coordinate of the system  we obtain the following parameters:
(note that we have normalized the singularities via the M\"obius transformation $z\mapsto z \cdot t_3$)
\begin{eqnarray*}
 (t_1,t_2,t_3)= (t,0,1),  & (\th_1,\th_2,\th_3,\th_4)=(c-\frac{1}{2}, a+c-1, c-\frac{1}{2}, a+c-1),
  \end{eqnarray*}
\begin{eqnarray*}\lambda =-b,\quad \tparam=b^2,& \Mu= \frac{-a-2c+2}{2b},&
 \nue=-\frac{3}{4b}. 
\end{eqnarray*}
The system (\ref{Schles}) 
using \eqref{Qi} reads:
{\small \begin{eqnarray*}
 q^1_{11}&=&\frac{( ba+a+2 b c-3 b-1)(a+2 c-2) }{4b(a+c-2)}, 
\\
 q^1_{12}&=&\frac{1}{b(b-1)}, \\
 q^1_{21}&=&\frac{(b-1)(ba+a-2 b+2 c-2) (-ba-a-2 b c+3 b+1)(a+2 c-2)(a-1)}{16b(a+c-2)^2}, \\
 q^1_{22}&=& -q^1_{11}+(c-\frac{1}{2}) \\
 q^2_{11}&=&\frac{(-b^2+2b-1)(a+2 c-2)(a-1)}{ 4b (a+c-2)},\\
 q^2_{12}&=&\frac{1}{b},\\
\end{eqnarray*}}
{\small \begin{eqnarray*}
 q^2_{21}&=&\frac{(b-1)^2((2-a-2 c)(a-1) (b^2+1)+(-2 a^2-4 a c+6 a-4 c^2+8 c-4) b)(a+2 c-2)(a-1)}{16b(a+c-2)^2},\\
 q^2_{22}&=&-q^2_{11}+(a-c-1), \\
q^3_{11}&=&\frac{(ba+a-b+2 c-3)(a+2 c-2)}{4(a+c-2)},\\
  q^3_{12}&=&-\frac{1}{b-1},\\
 q^3_{21} &=&\frac{(b-1)(ba+a-b+2 c-3)(ba+a+2 b c-2 b-2)(a+2 c-2)(a-1) }{16(a+c-2)^2},\\
 q^3_{22}&=&-q^3_{11}+(c-\frac{1}{2})\end{eqnarray*}
}
and
\begin{eqnarray*} -(Q_1+Q_2+Q_3)=\left( \begin{array}{cc}
                            -(a+2c-2) & 0 \\
                            0      & -c 
                    \end{array}\right).\end{eqnarray*}
The apparent singularity of the second coordinate $y_2$ is the zero
of $Q_{21}$:
\begin{eqnarray*} \nlambda&=&-b+\frac{(-2 a-2 c+4)(b+1)^2b}{(a-2)(a+2c-3)(b+1)^2+2 (a^2+2 ca-5 a+2 c^2-6 c+6)b}.
\end{eqnarray*}

\section{Middle convolution}
\label{middleconvolution}
For the convenience of the reader we give here 
a short review of the middle convolution for Fuchsian systems
(see  \cite{DR}).

Let $f(z)=( f_1(z),\ldots,f_n(z))^{\tr}$
be a solution of the Fuchsian system
\begin{eqnarray*} DY=AY=\sum_{i=1}^r \frac{A_i}{z-t_i} Y, \quad A_i  \in \Mat_n(\CC), \end{eqnarray*}
and
 $[t_i,y]=\gamma_{t_i}^{-1} \gamma_y^{-1}\gamma_{t_i}\gamma_{y}$ be a Pochhammer cycle around  $t_i$ and $y.$ 
Then the Euler transform of $f(z)$ with respect to $\mu\in\C$ and $[t_i,y]$
\begin{eqnarray*}  \left(\begin{array}{cc}
             \int_{[t_i,y]} f(z)(y-z)^{\mu} \frac{dz}{z-t_1} \\
                              \vdots \\ 
               \int_{[t_i,y]} f(z)(y-z)^{\mu} \frac{dz}{z-t_r} 
       \end{array}
                 \right),& i=1,\ldots,r, \end{eqnarray*}
is a solution of the (Okubo-) system
\begin{eqnarray*} (yI_{nr}-T) DY= BY:=(\left( \begin{array}{cccc}
                    A_1 & \ldots & A_r \\
                     & \vdots \\
                     A_1 & \ldots & A_r \\
                 \end{array}
                 \right) +\mu \;I_{nr}) Y,\quad T=\diag(t_1I_n,\ldots,t_rI_n)\end{eqnarray*}
\begin{eqnarray*}          \Leftrightarrow
         DY= \sum_{i=1}^r \frac{B_i}{y-t_i}Y,& B_i \in \Mat_{nr}(\CC), \end{eqnarray*}
where $I_n$ is the identity $n\times n$ matrix. 
In general this system is not irreducible and has the following two $\langle
B_1,\ldots,B_r \rangle$-invariant subspaces:
\begin{eqnarray*} \k =\oplus_{i=1}^r \ker(A_i), & \l= \ker(B)=\langle (v,\ldots,v)^\tr \mid v \in \ker(\sum_{i=1}^r A_i+\mu I_n)\rangle.\end{eqnarray*} 
Factoring out this subspace we obtain a Fuchsian system in dimension $m$
\begin{eqnarray*}  m=nr-\sum_{i=1}^r \dim(\ker A_i)-\dim(\ker(\sum_{i=1}^r A_i+\mu I_n)).
\end{eqnarray*}
Let $M_i$ be the monodromy of the system $DY=AY$ at $t_i.$
Then, if the system is irreducible and 
\begin{eqnarray*}
    \rk(A_i)&=&\rk(M_i-I_n),\quad i=1,\ldots,r \\
        \rk(\sum_{i=1}^r A_i+\mu I_n)&=&
     \rk(M_1\cdots M_r \lambda -I_n),\quad \lambda=e^{2  \pi i \mu} \end{eqnarray*}
then the factor system is again irreducible and it is called the { \it middle
convolution} of $DY=AY$ with $\mu.$ Thus if we start with a two dimensional system in Schlesinger form (\ref{Schles}),
where the parameters are $(\th_1,\ldots,\th_4)$ 
then there are in general two possibilities to get again a new  two dimensional system (\ref{Schles})
via the middle convolution.
We can apply the  middle convolution either with $\mu=\alpha$ or with $\mu =\alpha+\th_4-1$ and diagonalize the residue matrix at $\infty$ to
\begin{eqnarray*} \left( \begin{array}{cc}
                            \tilde{\alpha} & 0 \\
                            0      &  \tilde{\alpha}+\tilde{\th}_4-1
                    \end{array}\right).\end{eqnarray*}
This changes the parameters as follows:
\begin{eqnarray*} (\tilde{\th}_1,\ldots,\tilde{\th}_4)= (\th_1+\mu,\ldots,\th_3+\mu,\th_4-\mu+2\alpha),
&\tilde{\alpha}=-\mu. \end{eqnarray*}
Note that it is possible to write down also the transformation for
 the apparent singularity (see \cite[Section 5]{FiHa})
 which is also known as an Okamoto transformation.

\section{Application of the middle convolution I}
\label{midconvI}
In this section we show that the algebraic solution in  row 2 of Table 1 can be obtained via the middle convolution of the Schlesinger system corresponding to the algebraic solution in row 3 of Table 1.
 
We continue with our example by 
determining the middle convolution of $DY=QY$ with $\mu=-c$. 
Hence,  for the (Okubo-) system
$(zI_6-T) DY= BY \iff  DY= \sum_{i=1}^3 \frac{B_i}{y-t_i}Y,\;
 B_i \in Mat_6(\CC),$ we get 
    \begin{eqnarray*}B=\left( \begin{array}{ccc}
        Q_1 -c\cdot I_2& Q_2 & Q_3 \\
         Q_1 & Q_2-c\cdot I_2 & Q_3 \\
        Q_1 & Q_2 & Q_3 -c\cdot I_2
     \end{array} \right) \end{eqnarray*}
and $T=\diag (t,t,0,0,1,1).$
The $\langle B_1,B_2 ,B_3 \rangle $
invariant subspace has dimension 4, since
 \begin{eqnarray*}\ker Q_i=\langle \left( \begin{array}{c} 
              q^i_{1,2}         \\
             -q^i_{1,1} 
             \end{array}\right) \rangle &&
\l=\langle(0,1,0,1,0,1)^{\tr} \rangle.  \end{eqnarray*}
In order to get a 2 dimensional factor system we set
\begin{eqnarray*}
      S=\left(  \begin{array}{cccccc} 
            q^1_{1,2} &0 &0 &0 &0 & 0        \\
             -q^1_{1,1} &0 &0&1 &0 & 0     \\
            0 &  q^2_{1,2} &0 &0&1& 0      \\
            0 &-q^2_{1,1}  &0 &1 &0& 0     \\
            0 &0& q^3_{1,2} &0 &0&1       \\
            0 &0&- q^3_{1,1}& 1 &0& 0       
 \end{array} \right).
\end{eqnarray*}
Thus
\begin{eqnarray*} D(SY) = \sum_{i=1}^3 \frac{S B_i S^{-1}}{z-t_i}  (SY) =
        \left( \begin{array}{ccccccc} 
             \ast & 0 & 0 & 0 &\ast         \\
             0  &\ast & 0 &  0& \ast   \\
             0 &  0 &\ast & 0 &\ast      \\
             0 &  0 & 0  & 0 &  \ast  \\
            0  & 0& 0&  0 &\tilde{A}        
 \end{array}\right) (SY)
\end{eqnarray*}
and we get  the 2 dimensional factor system $DY=\tilde{A}Y=\sum_{i=1}^3 \frac{\tilde{A}_i}{z-t_i}Y$.
Transforming this system into Schlesinger form via
\[ Y\mapsto \tilde{S}Y,\quad 
 \tilde{S}=\left(  \begin{array}{cc}
         -(\sum_{i=1}^3 \tilde{A_i}+c\cdot I_2)_{1,2}&-(\sum_{i=1}^3 \tilde{A_i}-c\cdot I_2-a\cdot I_2+2\cdot I_2)_{1,2}\\
  
                (\sum_{i=1}^3 \tilde{A_i}+c\cdot I_2)_{1,1}&(\sum_{i=1}^3 \tilde{A_i}-c\cdot I_2-a\cdot I_2+2\cdot I_2)_{1,1}
 \end{array}\right),
\] 
where the columns of $\tilde{S}$ consist of eigenvectors of $\sum_{i=1}^3\tilde{A}_i$ with respect to the eigenvalues $-c$ and $c+a-2$,
   we obtain finally
\begin{eqnarray}
\label{dihedral} 
DY &=&A Y, \ A=\frac{A_1}{z-b^2}+\frac{A_2}{z}+\frac{A_3}{z-1}
\end{eqnarray}
\begin{eqnarray*} A_1&=&\frac{1}{4 b}\left(\begin{array}{cc}
     -a b-a-2 b c+b+1 &
      a b+a+2 b c-3 b-1\\
       -a b-a-2 b c+b+1&
       a b+a+2 b c-3 b-1
\end{array}\right),\\
 A_2&=&\frac{1}{4 b} \left(\begin{array}{cc}
   a b^2+2 a b+a-b^2-2 b-1 &
  a b^2-a-b^2+1\\
  -a b^2+a+b^2-1&
-a b^2+2 a b-a+b^2-2 b+1
\end{array}\right),\\
 A_3&=&\frac{1}{4}\left(\begin{array}{cc}
      -a b-a+b-2 c+1&
 -a b-a+b-2 c+3\\
 a b+a-b+2 c-1&
 a b+a-b+2 c-3
\end{array}\right)\end{eqnarray*}
and  
\begin{eqnarray*} -(A_1+A_2+A_3)&=&\left(\begin{array}{cc}
c  &
          0\\
          0&
           -(a+c-2)
\end{array}\right).\end{eqnarray*}
Computing the entries of $A$ we get
\begin{eqnarray*}
A_{11}&=& \frac{-4 c z^2+((1-a+2c)(b^2+1)+2(1-a)b) z+(a-1)(b+1)^2 b}{4 z(z-1)(z-b^2)}, \\
A_{12}&=& \frac{(b^2-1)((a+2 c-3) z+b(a -1))}{4 z(z-1)(z-b^2)}, \\
A_{21}&=& \frac{(b^2-1)((1-a-2 c) z+(b-a b))}{4 z(z-1)(z-b^2)}, \\
A_{22}&=& \frac{4(a+c-2) z^2+((5-3a-2c)(b^2+1)+2(a-1)b) z+(1-a)(b-1)^2 b}{4 z(z-1)(z-b^2)}.
\end{eqnarray*}
Thus  the parameters are 
\[(\th_1,\ldots,\th_4)=(-1/2,a-1,-1/2,-(a+2 c-3)).\]
The apparent singularity $\la_1$ for the first, resp. $\la_2$ for the second
coordinate is
\begin{eqnarray*} \la_1=(\frac{-a+1}{a+2c-3})b, &&
 \la_2=(\frac{-a+1}{a+2c-1})b. \end{eqnarray*}
Since $t=b^2$ we obtain the relation
\begin{eqnarray*} \th_4^2 \;\la_1^2 = t\;\th_2^2, &&
 (2-\th_4)^2 \;\la_2^2 = t\; \th_2^2.\end{eqnarray*}
Computing the \SL-form for the first coordinate (after transforming $Y\mapsto (z-1)^{1/2}(z-t)^{1/2}Y$ which changes the parameters $(\th_1,\ldots,\th_4)$ (but not the apparent singularities) to
$(1/2,a-1,1/2,-(a+2 c-3))$ we get a simple formula for $\mu_1$
\begin{eqnarray*} \mu_1=\frac{-a-2c+3}{2b}.\end{eqnarray*}

\section{Monodromy}
\label{monodromy}
In this section we show that (\ref{dihedral}) 
 arises also  as a pullback of hypergeometric
differential equations.
 By determining the monodromy group it turns
out that this monodromy group is finite and therefore
 by a well known result of Klein
the claim follows. 
Further we indicate that all the other families of Picard-Fuchs equations
are related to those with finite monodromy.

\begin{prop}
The projective monodromy group of \eqref{dihedral}
is (up to conjugation) contained in the orthogonal group $\GO_2(\CC)$.
Moreover, if the parameters $a$ and $c$ are rational numbers, i.e. if
\eqref{dihedral} is a Picard-Fuchs equation, then the monodromy group
is even finite, i.e. a dihedral group.
\end{prop}

\begin{proof}
To prove this statement let
$ M_t, M_0, M_1$
denote the monodromy at $t, 0 $ and $1.$ 
If the parameter $b$ tends to $1$  we see that the system becomes
reducible and monodromy group is abelian:

\begin{eqnarray*}A \to  \left(\begin{array}{cc}
\frac{-a-c z+1}{z^2-z} &
0\\
0&
\frac{a+c-2}{z-1}
\end{array}\right).\end{eqnarray*}
Hence $M_t M_1$ and $M_0$ commute and are diagonal matrices.
Since $M_t$ and $M_1$ are reflections, the group generated by them
is an orthogonal group. 
Hence $M_t$ and $M_1$ normalize also $M_0$ and the
claim follows.
\end{proof}

Applying the M\"obius transformation $z \mapsto \frac{1}{z}$ that
permutes the residue matrices $Q_i$ of the Schlesinger system corresponding to row 5 (resp. row 6), column 2 in Table~{1} and scaling the new $Q_2$ 
we get the following pairs of eigenvalues  for $Q_1, Q_2, Q_3,-(Q_1+Q_2+Q_3)$  
 \begin{eqnarray*}
 (a-\frac{1}{2},0),\;(a-\frac{3}{2},0),\;(a-\frac{1}{2},0),\; (-\frac{3a-3}{2},-\frac{3a-2}{2})
\end{eqnarray*}
resp.
\begin{eqnarray*}
 (a-\frac{1}{2},0),\;(-2a+2,0),\;(a-\frac{1}{2},0),\; (-\frac{1}{3},-\frac{2}{3}).
\end{eqnarray*}
The middle convolution with $\mu=-(a-1)$  yields a three dimensional Fuchsian system with the following triples of eigenvalues of the residue matrices
\begin{eqnarray*}
 (\frac{1}{2},0,0),\;(-\frac{1}{2},0,0),\;(\frac{1}{2},0,0),\;(-\frac{a-1}{2},-\frac{a}{2},a-1)
\end{eqnarray*}
resp.
\begin{eqnarray*}
(\frac{1}{2},0,0),\;(-3a+3,0,0),\;(\frac{1}{2},0,0),\; (-\frac{1}{3}+a-1,-\frac{2}{3}+a-1,a-1).
\end{eqnarray*}

Using the explicit construction for the middle convolution and similar arguments as in the above Proposition, one easily sees that the monodromy groups of these third order differential equations are finite imprimitive reflection groups contained
in $T\rtimes S_3,$ where $T$ denotes the group of diagonal matrices.
In the next section we show that in Table 1 row 4 and row 5 are also related via
the middle convolution. Hence all the families  of Picard-Fuchs equations
corresponding to Table~1  are related to those with finite monodromy. 
Since the corresponding Picard-Fuchs differential equation of row 1 in Table~{1} is the hypergeometric one it is well known that it is obtained via the convolution  of a one dimensional differential equation with finite monodromy.

\section{Application of the middle convolution II}
\label{midconvII}

As in the previous example we show that
the middle convolution relates the results in Table 1 row 4 and row 5.

The system (\ref{Schles}) in Schlesinger form corresponding to the
differential equation satisfied by $\int_\delta \frac{\dx}{y}$, where $y$ is
from row 4 in Table 1, reads:

{\small \begin{eqnarray*}
q^1_{1,1}&=&\frac{(18 a^2 b^2+18 a^2 b+18 a^2-3 a b^3-30 a b^2-48 a b-36 a+2 b^3+12 b^2+24 b+16)}{18 a b+18 a-27 b-27} \\
   q^1_{1,2}&=&\frac{-b^4-2 b^3}{b^2-1}\\
 q^1_{2,1}&=& \frac{q^1_{11} q^1_{22}}{q^1_{12}}\\
   q^1_{2,2}&=&(a-\frac{1}{2})-q^1_{1,1}\\
   q^2_{1,1}&=&\frac{(3a-2)(b^2+(-6a+7)b+1)(b-1)^2}{18 a b^2-27 b^2} \\
   q^2_{1,2}&=&b^2+2 b\\
q^2_{2,1}&=& \frac{q^2_{11} q^2_{22}}{q^2_{12}}\\
   q^2_{2,2}&=&(3a-\frac{3}{2})-q^2_{1,1}\\
   q^3_{1,1}&=&\frac{(18 a^2 b^3+18 a^2 b^2+18 a^2 b-36 a b^3-48 a b^2-30 a b-3 a+16 b^3+24 b^2+12 b+2)}{18 a b^3+18 a b^2-27 b^3-27 b^2} \\
   q^3_{1,2}&=&\frac{b^2+2 b}{b^2-1}\\
   q^3_{2,1}&=&\frac{q^3_{11} q^3_{22}}{q^3_{12}}\\
   q^3_{2,2}&=&(a-\frac{1}{2})-q^3_{1,1}
\end{eqnarray*}}
\begin{eqnarray*}
-(Q_1+Q_2+Q_3)&=&\left(\begin{array}{cc}
  -3 a+2&
        0\\
         0&
 -2 a+\frac{1}{2}
\end{array}\right).
\end{eqnarray*}

Applying the middle convolution to $DY=QY$ with $\mu=-(3a-2)$
 and transforming the 2 dimensional factor system into Schlesinger form with singularities at $t, 0, 1$ and $\infty$ we get
the system (\ref{Schles}) with 
\begin{eqnarray*}
q^1_{1,1}&=&\frac{(3 a-2)(b+(-6 a+2))(b+2)^2}{ 9 (4a-1) (b+1)}
\\
q^1_{1,2}&=&\frac{b+2}{(72 a-18)(b+1)}
\\
q^1_{2,1}&=&\frac{(3 a-2)((6 a-4) b^3+(-36 a^2+60 a-24) b^2+(24 a-21) b-5)(-b+(6 a-2))(b+2)}{9 (4a-1) (b+1)}
\\
q^1_{2,2}&=&(-2 a+\frac{3}{2})-q^1_{1,1}\\
q^2_{1,1}&=&\frac{(-3 a+2) (b^2+(-6 a+4) b+1) (b^2+b+1)}{9 (4a-1) b^2}
\\
q^2_{1,2}&=&\frac{-(b^2+b+1) }{ 18 (4a-1) b^2}
\\
q^2_{2,1}&=&\frac{( (4-6 a) b^4+(36 a^2-54 a+20) b^3+(36 a^2-96 a+33) b^2+(36 a^2-54 a+20) b+4-6 a)}{9 (4a-1) b^2 } \cdot \\
 &&(3 a-2)(-b^2+ (6 a-4) b-1)
\\
q^2_{2,2}&=&\frac{1}{2}-q^2_{1,1}\\
q^3_{1,1}&=&\frac{(3 a-2)((-6 a+2) b+1)(2 b+1)^2}{9(4a-1)b^2(b+1)}
\\
q^3_{1,2}&=&\frac{2b+1}{18(4a-1)b^2(b+1)}
\\
q^3_{2,1}&=&\frac{(3 a-2)(-5 b^3+(24 a-21) b^2+(-36 a^2+60 a-24) b+(6 a-4))((6 a-2) b-1)(2 b+1)}{9(4a-1)b^2(b+1)}
\\
q^3_{2,2}&=&(-2 a+\frac{3}{2})-q^3_{1,1}
\end{eqnarray*}

\begin{eqnarray*}
 -(Q_1+Q_2+Q_3)=\left(\begin{array}{cc}
          3a-2&0\\
                 0&  a-\frac{3}{2}
                  \end{array}\right)=
  \left(\begin{array}{cc}
          3a-2&0\\
                 0& 3a-2+(-2 a+\frac{3}{2})-1
                  \end{array}\right).     
\end{eqnarray*}

We obtain  the parameters $(\th_1,\ldots,\th_4)=
(-2 a+\frac{3}{2},\frac{1}{2},-2 a+\frac{3}{2},-2 a+\frac{3}{2})$.
The apparent singularity for the first coordinate is
\begin{eqnarray*} \lambda&=&\frac{b^2+b+1}{b^3+2 b^2}.\end{eqnarray*}
Using that 
$t=\frac{2 b+1}{b^4+2 b^3}$ we see that $\lambda$ and $t$ satisfy
\begin{eqnarray*} \lambda ^4-2 t \lambda ^3-2 \lambda ^3+6 t \lambda ^2-2 t^2 \lambda -2 t \lambda +t^3-t^2+t=0.
\end{eqnarray*}
Since this is also the relation  for
\begin{eqnarray*} t=\frac{ b^4-6 b^2-8 b-3}{b^4-6 b^2+8 b-3}, &&
                \lambda=\frac{b^3+b^2+3 b+3}{b^3+b^2-5 b+3}
\end{eqnarray*}
from Table 1 row 5 the claim follows.

{}


\begin{thebibliography}{1}

\bibitem{ben05}
B. Ben~Hamed and L. Gavrilov.
\newblock Families of {P}ainlev\'e {VI} equations having a common solution.
\newblock {\em Int. Math. Res. Not.}, (60):3727--3752, 2005.


\bibitem{Boalch03}
P. Boalch.
\newblock Painlevé equations and complex reflections. 
\newblock {\em Ann. Inst. Fourier} 53 no. 4, 1009--1022, 2003.


\bibitem{boa06}
P. Boalch.
\newblock The fifty-two icosahedral solutions to {P}ainlev\'e {VI}.
\newblock {\em J. Reine Angew. Math.}, 596:183--214, 2006.


\bibitem{CL07}
S. Cantat and F. Loray.
\newblock Holomorphic dynamics, Painlevé VI equation and Character Varieties.
\newblock {\em arXiv :0711.1579v2}, 1-69, 2007.  


\bibitem{DR}
M. Dettweiler and S. Reiter.
\newblock  Middle convolution of Fuchsian systems and the construction of rigid differential systems.  
\newblock {\em J. Algebra}  318,  no. 1, 1--24, 2007.


\bibitem{dor01}
C.~F. Doran.
\newblock Algebraic and geometric isomonodromic deformations.
\newblock {\em J. Differential Geom.}, 59(1):33--85, 2001.

\bibitem{FiHa}
Y. Haraoka and  G. Filipuk,
\newblock Middle convolution and deformation for Fuchsian systems.  
\newblock {\em J. Lond. Math. Soc. (2)}  76,  no. 2, 438--450, 2007.



\bibitem{GPS01}
G.-M. Greuel, G.~Pfister, and H.~Sch\"onemann.
\newblock {\sc Singular} 2.0.
\newblock {A Computer Algebra System for Polynomial Computations}, Centre for
  Computer Algebra, University of Kaiserslautern, 2001.
\newblock {\tt http://www.singular.uni-kl.de}.


\bibitem{gri66}
P.~A. Griffiths.
\newblock The residue calculus and some transcendental results in algebraic
  geometry. {I}, {II}.
\newblock {\em Proc. Nat. Acad. Sci. U.S.A.}, 55:1303--1309; 1392--1395, 1966.




\bibitem{her91}
S. Herfurtner.
\newblock Elliptic surfaces with four singular fibres.
\newblock {\em Math. Ann.}, 291(2):319--342, 1991.

\bibitem{iwa91}
K. Iwasaki, H. Kimura, S. Shimomura, and M. Yoshida.
\newblock {\em From {G}auss to {P}ainlev\'e}.
\newblock Aspects of Mathematics, E16. Friedr. Vieweg \& Sohn, Braunschweig,
  1991.
\newblock A modern theory of special functions.

\bibitem{kit05}
A.~V. Kitaev.
\newblock Grothendieck's dessins d'enfants, their deformations, and algebraic
  solutions of the sixth {P}ainlev\'e and {G}auss hypergeometric equations.
\newblock {\em Algebra i Analiz}, 17(1):224--275, 2005.

\bibitem{mov08}
H.~Movasati.
\newblock On {R}amanujan relations between {E}isenstein series.
\newblock {\em Submitted}, 2008.

\bibitem{shi04}
H. Shiga, T. Tsutsui, and J. Wolfart.
\newblock Triangle {F}uchsian differential equations with apparent
  singularities.
\newblock {\em Osaka J. Math.}, 41(3):625--658, 2004.
\newblock With an appendix by Paula B. Cohen.











\end{thebibliography}

\end{document}